\documentclass[number,citesort,seceqn,dvips]{arxbj}
\usepackage{graphicx}

\aid{0}
\volume{17}
\issue{4}
\pubyear{2011}
\firstpage{1368}
\lastpage{1385}
\doi{10.3150/10-BEJ320}

\makeatletter

\newtheorem{corollary}[prop]{Corollary}
\newtheorem{lemma}{Lemma}[section]
\newtheorem{prop}{Proposition}[section]
\newproclaim{definition}[prop]{Definition}
\newremark{remark}[lemma]{Remark}
\newtheorem{theorem}[prop]{Theorem}

\newcommand{\wh}{\hat}
\newcommand{\eps}{\varepsilon}
\newcommand{\RR}{\mathbb R}
\newcommand{\EE}{\mathbb{E}}
\newcommand{\PP}{\mathbb{P}}

\newcommand{\Xcal}{\mathcal X}
\newcommand{\argmin}{\operatorname{arg\,min}}
\makeatother

\begin{document}
\begin{frontmatter}

\title{Support vector machines with a reject option}
\runtitle{Reject SVM}

\begin{aug}
\author[1]{\fnms{Marten} \snm{Wegkamp}\thanksref{1}\ead[label=e1,text={marten.wegkamp@cornell.edu}]{marten.wegkamp@cornell.edu}} \and
\author[2]{\fnms{Ming} \snm{Yuan}\thanksref{2}\corref{}\ead[label=e2]{myuan@isye.gatech.edu}}
\runauthor{M. Wegkamp and M. Yuan}
\address[1]{Department of Mathematics and Department of Statistical
Science,
Cornell University, Ithaca, NY 14853, USA. \printead{e1}}
\address[2]{Milton Stewart School of Industrial and Systems Engineering,
Georgia Institute of Technology, Atlanta, GA 30332, USA. \printead{e2}}
\end{aug}

\received{\smonth{6} \syear{2009}}
\revised{\smonth{1} \syear{}2010}

%
\begin{abstract}
This paper studies $\ell_1$ regularization with high-dimensional
features for support vector machines with a~built-in reject option
(meaning that the decision of classifying an observation can be
withheld at a cost lower than that of misclassification). The procedure
can be conveniently implemented as a linear program and computed using
standard software. We prove that the minimizer of the penalized
population risk favors sparse solutions and show that the behavior of
the empirical risk minimizer mimics that of the population risk
minimizer. We also introduce a notion of classification complexity and
prove that our minimizers adapt to the unknown complexity. Using a
novel oracle inequality for the excess risk, we identify situations
where fast rates of convergence occur.
\end{abstract}

%
\begin{keyword}
\kwd{adaptive prediction}
\kwd{classification with a reject option}
\kwd{lasso}
\kwd{oracle inequalities}
\kwd{sparsity}
\kwd{support vector machines}
\kwd{statistical learning}
\end{keyword}

\end{frontmatter}

\section{Introduction}

In this paper we further investigate the new classification rules
introduced in \cite{BW08,W07} with a built-in reject option in the
standard binary classification setting, where we observe independent
realizations $(X_i,Y_i)$, $i=1,\ldots, n$, of a random pair $(X,Y)$ in
$\mathcal X \times\{-1,+1\}$ (here, $\mathcal X$ is an arbitrary
space). A discriminant function $f\dvtx\mathcal X\to\RR$ classifies an
observation $x\in\Xcal$ into one of two classes, labeled $-1$ or $+1$.
Viewing $f(x)$ as a proxy value of the conditional probability $
\eta(x) = \PP\{ Y=1| X=x\}, $ we are less confident for small values of
$|f(x)|$, corresponding to $\eta(x)$ near 1$/$2. Our strategy is to
report $\operatorname{sgn}(f(x))\in\{-1,1\}$ if $|f(x)|$ exceeds some
prescribed threshold $\tau$ and withhold decision otherwise. Assuming
that the cost of making a wrong decision is 1 and that of withholding a
decision is $d$, the appropriate risk function is
\[
R_\ell(f) = \EE[ \ell(Yf(X)) ] = \PP\{ Y f(X) < -\tau\} + d \PP\{ | Y
f(X) | \le\tau\}\label{risk_orig}
\]
with the discontinuous loss function
\[
\ell(z) = \cases{
1, &\quad if $z<-\tau$,\cr
d,&\quad if $|z|\le\tau$,\cr
0, &\quad otherwise. }\label{dc}
\]
Since we always reject if $d=0$ and never reject if $d \ge1/2$ (see
\cite{HW06}), we take $0< d\le1/2$ in what follows without loss of
generality. Although the minimizer of this risk is not unique, all such
minimizers correspond to the unique classification rule that assigns
$-1,+1$ or withhold decision, depending on which of $1-\eta$, $\eta$ or
$d$ is smallest. The smallest risk is $\mathbb E[
\min\{\eta(X),1-\eta(X),d\}]$ and we may interpret the cost $d$ as the
largest conditional probability of misclassification that is considered
tolerable.

In practice, minimization of the empirical counterpart $\wh R_\ell(f)=
(1/n) \sum_{i=1}^n \ell(Y_i f(X_i))$ of $R_\ell(f)$ over a large class
of functions $f$ is computationally not feasible. For this reason, we
could replace the loss function $\ell$ by a convex surrogate loss
function and consider discriminant functions $f$ of the form $ {\sf
f}_\lambda(x)= \sum_{j=1}^M \lambda_j f_j(x) $ based on a set of known
functions $f_j\dvtx\mathcal X\to\RR$ and coefficients $\lambda_j\in\RR$,
$1\le j\le M$. Following \cite{BW08}, we will consider the
generalized hinge loss
\[
\phi(z) =
\cases{
1-a z, &\quad if $z<0$,\cr
1- z, &\quad if $0\le z<1$,\cr
0, &\quad otherwise
}\label{hinge}
\]
with slope $a=(1-d) / d > 1$. Observe that
$\phi(z)$ is piecewise linear, so that minimization of the empirical risk
%
%
\begin{equation}\label{emp}
 \wh R_\phi({\sf f}_\lambda)= \frac1n \sum_{i=1}^n \phi(
Y_i {\sf f}_\lambda(X_i))
\end{equation}
can be solved by a tractable linear program. Crucial for the choice of
$\phi(z)$ is that it is classification calibrated: the unique minimizer
\[
f_0(x) = \cases{
-1, &\quad  if $\eta(x) < d,$\cr
 0,&\quad  if $d\le\eta(x)\le1-d,$\cr
+1, &\quad  if $\eta(x) > 1-d$
}\label{bayes_d}
\]
of $R_\phi(f)=\EE[ \phi(Y f(X)) ]$ also minimizes the risk $R_\ell(f)=
\EE[ \ell(Y f(X)) ]$ over all measurable $f\dvtx\mathcal X\to\RR$ for all
$\tau<1$; see, for example, \cite{BW08,YW09}.

At this point it is important to note that truncating the minimizer
$\operatorname{sgn}(2\eta-1)$ of the hinge-loss-based risk $\EE(
1-Yf(X) )_+$ does not yield the optimal rule for any positive threshold
$\tau$. This is the reason why we generalize the hinge loss instead. In
addition to the generalized hinge loss, there are also other choices of
the surrogate loss function and corresponding truncation value $\tau$
that are classification calibrated. The treatment for the generalized
hinge loss differs considerably from that for other losses, such as the
logistic, exponential and quadratic loss, which are smoother. We refer
to \cite{YW09} for a detailed discussion.

Observe that $\phi(z)\ge\ell(z)$ for all $\tau\le1-d$ and,
subsequently, $\EE[ \ell(Y f(X)) ] \le\EE[ \phi(Y f(X)) ]$. It is shown
in \cite{BW08} that a similar relationship remains true for the excess
risks, that is, the inequality
\[
\EE[ \ell(Y f(X)) ] - \EE[ \ell(Y f_0(X)) ] \le \EE[ \phi(Y f(X)) ] -
\EE[ \phi(Y f_0(X)) ]
\]
holds for all $d\le\tau\le1-d$. This property is useful for deriving
oracle inequalities in terms of the $\ell$-risk since minimization of
(\ref{emp}) produces oracle inequalities in terms of the $\phi$-risk
rather than the $\ell$-risk directly.

Of particular interest here is the case where the number of basis
functions, $M$, is large when compared with the sample size $n$.
Usually, the minimization of the empirical risk $\wh R_\phi({\sf
f}_\lambda)$ is computed under a restriction on the quadratic term $
\sum_{j=1}^M \lambda_j^2$. Here, we opt instead for an $\ell_1$-type
restriction $\|\lambda\|_{\ell_1}:=\sum_{j=1}^M |\lambda_j| $ and
estimate $\sf f_{\lambda}$ by $\sf f_{\wh{\lambda}(r)}$, where
%
%
\begin{equation}
\wh{\lambda}(r):=\argmin\limits_{\lambda\in\RR^M}\bigl(\wh R_\phi({\sf f}_\lambda) +
r\|\lambda\|_{\ell_1}\bigr)
\end{equation}
and $r>0$ is a tuning parameter. The choice of an $\ell_1$ penalty
reflects our preference for sparse solutions, which is desirable when
$M$ is large.

In the remainder of this paper, we study the properties of
$\wh{\lambda}(r)$ and its population counterpart,
%
%
\begin{equation}
\lambda(r):=\argmin\limits_{\lambda\in\RR^M}\bigl(R_\phi({\sf f}_\lambda) +
r\|\lambda\|_{\ell_1}\bigr).
\end{equation}
We establish oracle inequalities for $\lambda(r)$ and $\wh{\lambda}(r)$
in Sections \ref{sec2} and \ref{sec3}, respectively. The results that we obtain are
similar in spirit to those from \cite{K09,TG06,W07}. However,
\cite{TG06,W07} do not discuss properties of $\lambda(r)$, and our
results in Section \ref{sec2} obtained here extend those proved by \cite{K09} in
the context of twice differentiable loss functions. Furthermore, the
oracle inequalities for the penalized empirical risk minimizer
$\wh{\lambda}(r)$ in Section \ref{sec3} are much sharper than earlier results
from~\cite{W07} for $0\le d\le1/2$ and \cite{TG06} for $d=1/2$.
In particular, the new inequality reveals that the rate of convergence
of the excess risk of ${\sf f}_{\wh{\lambda}}$ can be even faster than
$1/n$ if the optimal discriminant function $f_0$ can be written as a
linear combination of the $f_j$'s in the dictionary. Moreover, we relax
the condition on the dictionary and do not require that the parameter
$\lambda$ is bounded. We emphasize that our results hold, in
particular, for $d=1/2$, the case of support vector machines without a~reject option, and generalize and extend the results obtained in \cite{TG06}.
In addition, novel empirical bounds on the error and reject
rate are given.
To demonstrate the feasibility of the $\ell_1$-regularized support
vector machine with a reject option, in Section \ref{sec4} we formulate $\wh
\lambda(r)$ as a solution of a~linear program and report some numerical
experiments. Some technical lemmas and a maximal inequality for a
weighted empirical process are collected in the \hyperref[appendix]{Appendix}.

\section{Properties of the theoretical solution}\label{sec2}

We begin by studying $\lambda(r)$, the population version of
$\wh{\lambda}(r)$. Recall that $\lambda(r)$ is defined by
%
%
\begin{equation}
\lambda(r) = \arg\min_{\lambda\in\RR^M} \{ R_\phi({\sf
f}_\lambda) + r
\|\lambda\|_{\ell_1}\}.
\end{equation}
In particular, $\lambda(0)$ minimizes the risk $ R_\phi({\sf
f}_\lambda)$ over ${\lambda\in\RR^M}$. By definition, we find that
%
%
\begin{equation}\label{def}
 R_\phi\bigl({\sf f}_{\lambda(r)} \bigr)+ r\|\lambda(r)\|_{\ell_1}
\le R_\phi({\sf f}_{\lambda}) + r\|\lambda\|_{\ell_1}
\end{equation}
holds for all $\lambda\in\RR^M$. This inequality applied to
$\lambda=\lambda(0)$ has the following consequences.

\begin{prop}\label{prop1}
Let $I_0=\{i\dvtx \lambda_i(0)\ne0\}$ be the support of $\lambda(0)$.
\begin{longlist}
\item[(a)] If $\| \lambda(0)\|_{\ell_1} = o(1/r)$ as $r\to0$, then
$R_\phi(\lambda(r))\to R_\phi(\lambda(0))$ as $r\to0$.
\item[(b)]
$\|\lambda(r)\|_{\ell_1}\le\|\lambda(0)\|_{\ell_1}$ for all $r>0$.
\item[(c)] $\sum_{j\notin I_0}| \lambda_j(r)-\lambda_j(0)|\le
\sum_{j\in I_0}|\lambda_j(r)-\lambda_j(0)|$.
\end{longlist}
\end{prop}
\begin{pf}
After applying inequality (\ref{def}) to $\lambda=\lambda(0)$ and
using the fact that $R_\phi({\sf f}_{\lambda(0)} )\le R_\phi({\sf
f}_{\lambda(r)})$, we get
\[
0\le R_\phi\bigl({\sf f}_{\lambda(r)} \bigr)-R_\phi\bigl({\sf f}_{\lambda(0)} \bigr)
\le r\|
\lambda(0)\|_{\ell_1}-r\|\lambda(r)\|_{\ell_1}\le r\|\lambda(0)\|
_{\ell_1},
\]
which implies (a). The second claim follows from
\[
R_\phi\bigl({\sf f}_{\lambda(r)}\bigr)+ r\|\lambda(r)\|_{\ell_1} \le R_\phi
\bigl({\sf
f}_{\lambda(0)}\bigr) + r\|\lambda(0)\|_{\ell_1}
\le R_\phi\bigl({\sf f}_{\lambda(r)}\bigr)+ r\|\lambda(0)\|_{\ell_1}.
\]
For the
proof of part (c), we first observe that $\|\lambda(r)\|_{\ell_1}\le
\|\lambda(0)\|_{\ell_1}$ is equivalent to
\[
\sum_{j\notin I_0} |\lambda_j(r) | \le\sum_{j\in I_0} |\lambda_j(0)|
- \sum_{j\in I_0} |\lambda_j(r)|.
\]
Next, we note that the term on the left equals $ \sum_{j\notin I_0} |
\lambda_j(0)-\lambda_j(r)|$ and we bound the term on the right by
$\sum_{j\in I_0} |\lambda_j(0)-\lambda_j(r) |$ using the triangle
inequality. This proves part (c).
\end{pf}

This result gives a simple condition for $R_\phi({\sf
f}_{\lambda(r)})\to R_\phi({\sf f}_{\lambda(0)})$ and shows that the
$\ell_1$ norm of the solution $\lambda(r)$ is always smaller than the
$\ell_1$ norm of $\lambda(0)$. Similar properties are established
by~\cite{K09} for minimizers of twice differentiable loss functions $\phi$
and $\ell_p$ norms for $p>1$. In contrast, we consider here a
non-differentiable loss function $\phi$ and $p=1$.

Our target is a sparse vector $\theta\in\RR^M$ with risk $R_\phi
({\sf
f}_\theta)$ close to $R_\phi({\sf f}_{\lambda(0)})$. Before we make
this precise, we need to introduce a few concepts depending on the
behavior of $\eta(X)$ near $d$ and $1-d$, and the set of functions
$f_j$.

\begin{definition}[(Classification complexity)]
The classification complexity is defined
as the largest number $\alpha\ge0$ such that, for some $A\ge1$ and all
$t>0$,
\[
\PP\{|\eta(X)-d|\le t\} \le A t^\alpha\quad \mbox{and}\quad \PP\{|\eta(X)-
(1-d)|\le t\}\le A t^\alpha.\label{sasha}
\]
\end{definition}

This notion of complexity is a generalization of Tsybakov's margin
condition \cite{T04} for $d=1/2$. The behavior of $\eta(X)$ is
obviously not relevant in the interval $(d,1-d)$, only at the
endpoints $ d$ and $1-d$. The inequality always holds for $\alpha=0$
and $A=1$. In contrast, $\alpha=+\infty$ describes the easiest
classification situation where we essentially require that $\eta(X)$
stays away from $d$ and $1-d$ with probability one. If $\eta(X)$ has a
density in the neighborhood of $d$ and $1-d$, then we have that
$\alpha=1$.

\begin{definition}[(Restricted eigenvalue condition)] Let $\theta\in\RR^M$, $c\ge1$ and
$\Psi$ be the $M\times M$ matrix with entries $ \Psi_{i,j} =4 \EE [
f_i(X) f_j(X) \omega(X) ]$ with $\omega(X)=\eta(X) \{ 1-\eta(X)\}$. For
$I=\{i\dvtx \theta_i\neq0\}$, the support of $\theta$, we define
\[
\kappa^2 (\theta,c)= \inf_{\lambda\ne\theta\in\RR^M\dvtx \|( \theta-
\lambda)_{I^C}\|_{\ell_1}\le c \|(\theta-\lambda)_{I}\|_{\ell_1} }
\frac { (\theta- \lambda)^\prime\Psi(\theta-\lambda) }{4\| (\theta-
\lambda )_I \|_{\ell_2}^2} .
\]
\end{definition}

The condition $\kappa(\theta,c)>0$ is a restrictive eigenvalue
condition on the Gram matrix $\Psi$ of the type introduced in
\cite{BRT} in the context of linear regression. Using similar
reasoning as in \cite{BRT}, page~1714, it is implied by the local
mutual coherence condition used in \cite{W07}.
We are now in position to state an oracle inequality for the excess
risk,
%
%
\begin{equation}
\Delta R_\phi\bigl({\sf f}_{\lambda(r)}\bigr):= R_\phi\bigl({\sf f}_{\lambda(r)}\bigr) -
R_\phi(f_0),
\end{equation}
of the regularized minimizer $\lambda(r)$ and the $\ell_1$-distance
between the vectors $\lambda(r)$ and $\theta$.

\begin{theorem}\label{theorie} Let $\alpha$ be the classification
complexity, and $\theta$ be such that $R({\sf f}_{\theta})\le R({\sf
f}_{\lambda(r)})$ and $\kappa=\kappa(\theta,1)>0$. Then, for any
%
%
\begin{equation}\label{reg2.4}
r\le(2C_F)^{-(2+\alpha)/\alpha} \{ 4A (2d)^\alpha \}
^{-1/\alpha} ( \kappa^{-2}\| \theta\|_{\ell_0} )^{
-(1+\alpha)/\alpha}
\end{equation}
with $C_F=\max_j\! \| f_j\|_{\infty}=\max_j \sup_x\! |f_j(x)|$ and
$\|\theta\|_{\ell_0}= \sum_{j=1}^M I\{ \theta_j\ne0\}$, we have
%
%
\begin{eqnarray}\label{theorem:oracle.th}
&&\Delta R_\phi\bigl({\sf f}_{\lambda(r)}\bigr) + r
\|\lambda(r)-\theta\|_{\ell_1}\nonumber\\ [-8pt]\\ [-8pt]
&&\quad\le 3 \Delta R_\phi({\sf f}_\theta)+ 6 \{ 4A (2d)^{\alpha} \}^{1/(2+\alpha)} \| {\sf f}_{\theta} -
f_0\|_{\infty} ( \kappa^{-2} r^2 \| \theta\| _{\ell _0} )
^{(1+\alpha)/(2+\alpha)}.\nonumber
\end{eqnarray}
\end{theorem}
\begin{pf}
Set $\delta=\lambda(r) -\theta$. Let $I=\{i\dvtx\theta_i\neq0\}$ be the
support of $\theta$. It is straightforward to derive from Proposition
\ref{prop1} that
\[
R_\phi\bigl({\sf f}_{\lambda(r)}\bigr)+ r\|\delta\|_{\ell_1} \le R_\phi({\sf
f}_{\theta}) + 2r \|\delta_{I}\|_{\ell_1}
\]
and, subsequently, that
\[
 r \| \delta_{ I ^C}
\|_{\ell_1} \le R_\phi({\sf f}_{\theta} )-R_\phi\bigl({\sf
f}_{\lambda(r)}\bigr) + r \| \delta_{I} \|_{\ell_1}\le
r\|\delta_{I}\|_{\ell_1}.
\]
The first inequality, combined with the assumption
$\kappa=\kappa(\theta,1)>0$, yields
\begin{eqnarray*}
\Delta R_\phi\bigl({\sf f}_{\lambda(r)}\bigr) + r \|\delta\|_{\ell_1} &\le&
\Delta R_\phi({\sf f}_{\theta}) + \kappa^{-1} \| {\sf f}_\delta\|
( r^2 | I| )^{1/2} \\
&\le& \Delta R_\phi({\sf f}_{\theta}) + \kappa^{-1} \| {\sf f}_\lambda
- f_0 \| ( r^2 | I| )^{1/2} + \kappa^{-1}\| {\sf f}_{\theta} - f_0 \| (
r^2 | I| )^{1/2},
\end{eqnarray*}
using the notation $ \| {\sf f}\| = \EE^{1/2} [ {\sf f}^2(X)\omega(X) ]
$ and $\omega(X)= \eta(X)(1-\eta)(X)$. By Lemma \ref{T} in
Appendix \ref{appendixA}, we find that
\[
 \| {\sf
f}_\lambda- f_0\| ^{2+2\alpha} \le 4A (2d)^{\alpha} \| {\sf
f}_\lambda- f_0\|_{\infty} ^{2+\alpha} \{ \Delta R_\phi({\sf
f}_{\lambda})\}^\alpha
\]
for $\lambda=\theta$ and $\lambda=\lambda(r)$. After we plug this bound
into the right-hand side of the previous display, we find
that
\begin{eqnarray*}
&&\Delta R_\phi\bigl({\sf f}_{\lambda(r)}\bigr) + r \|\delta\|_{\ell_1}\\
&&\quad\le \Delta R_\phi({\sf f}_{\theta}) + \kappa^{-1} ( r^2 | I| )^{1/2}
\{ 4A (2d)^\alpha\}^{1/(2+2\alpha)} \bigl\| {\sf f}_{\lambda(r)} - f_0
\bigr\| _{\infty} ^{(2+\alpha )/(2+2\alpha )} \bigl\{\Delta R_\phi\bigl({\sf
f}_{\lambda(r)}\bigr)\bigr\}^{\alpha
/(2+2\alpha)} \\
&&\qquad{}  + \kappa^{-1} ( r^2 | I| )^{1/2} \{ 4A
(2d)^\alpha\}^{1/(2+2\alpha)} \| {\sf f}_{\theta} - f_0 \|
_{\infty} ^{(2+\alpha)/(2+2\alpha)} \{\Delta R_\phi({\sf
f}_{\theta})\}^{\alpha/(2+2\alpha)}.
\end{eqnarray*}
Next, we apply Young's algebraic inequality,
\[
ab\le\frac{{a}^p}{p}+ \frac{b^q}{q}\qquad  \mbox{with }  p>1\quad \mbox{and}\quad q=\frac
{p}{p-1}\qquad \mbox{for all } a,b>0,
\]
to the last two terms on the right-hand side, with
$p=(2+2\alpha)/\alpha$ and $ q=(2+2\alpha)/(2+\alpha),$ to get
\begin{eqnarray*}
&&\Delta R_\phi\bigl({\sf f}_{\lambda(r)}\bigr) + r \|\delta\|_{\ell_1}\\
&&\quad\le \Delta R_\phi({\sf f}_{\theta}) + \frac{\alpha}{2+2\alpha} \bigl\{
\Delta R_\phi\bigl({\sf f}_{\lambda(r)}\bigr) + \Delta R_\phi({\sf f}_{\theta})
\bigr\} \\
&&\qquad{}+ \frac{2+\alpha}{2+2\alpha}\{ 4A
(2d)^\alpha\}^{1/(2+\alpha)}
(\kappa^{-2}r^2|I|)^{(1+\alpha)/(2+\alpha)}  \bigl(\bigl \| {\sf
f}_{\lambda(r)} - f_0\bigr\|_{\infty} +\| {\sf f}_{\theta} - f_0\|_{\infty}
\bigr).
\end{eqnarray*}
Since $\| {\sf f}_{\lambda(r)} - f_0 \|_{\infty} \le\| {\sf
f}_{\theta} - f_0\|_{\infty} + C_F \|\delta\|_{\ell_1},$ we deduce,
after invoking (\ref{reg2.4}), that
\begin{eqnarray*}
&&(2+\alpha)\Delta R_\phi\bigl({\sf f}_{\lambda(r)}\bigr) + (1+3\alpha/2)r \|
\delta\|_{\ell_1}\\
&&\quad\le (2+3\alpha)\Delta R_\phi({\sf f}_{\theta}) + 2(2+\alpha )\{ 4A
(2d)^\alpha\}^{1/(2+\alpha)}
(\kappa^{-2}r^2|I|)^{(1+\alpha)/(2+\alpha)} \| {\sf f}_{\theta} -
f_0\|_{\infty},
\end{eqnarray*}
and the conclusion follows.
\end{pf}

It is interesting to see that the bound (\ref{theorem:oracle.th})
crucially depends on the classification complexity parameter $\alpha$
and $\| {\sf f}_{\theta} - f_0\|_{\infty}$. In particular, if $f_0$
can itself be represented as a linear combination of the basis
functions, then $f_0 = {\sf f}_{\lambda(0)}$. In this case, provided
that $\kappa(\lambda(0),1)>0$, Theorem \ref{theorie} implies that
$\Delta R_\phi({\sf f}_{\lambda(r)}) + r
\|\lambda(r)-\lambda(0)\|_{\ell_1} \le0.$ In other words, we have the
following corollary.
\begin{corollary}
If $f_0 = {\sf f}_{\lambda(0)}$ and $\kappa(\lambda(0),1)>0$, then
$\lambda(r)=\lambda(0)$ for any
\begin{eqnarray*}
r\le(2C_F)^{-(2+\alpha)/\alpha} \{ 4A (2d)^\alpha \}
^{-1/\alpha} ( \kappa^{-2}\| \lambda(0)\|_{\ell_0} )^{
-(1+\alpha)/\alpha}.
\end{eqnarray*}
\end{corollary}

\section{$\ell_1$-regularized empirical generalized hinge risk
minimizers}\label{sec3}

In this section we study the estimate $\wh{\lambda}(2r)$. In what
follows, we will simplify notation so as not to show dependence of $\wh
\lambda$ on $r$ whenever no confusion occurs. Again, we emphasize that
our results hold, in particular, for $d=1/2$, the case of a support
vector machine without a reject option.

Note that the inequality
%
%
\begin{equation}\label{lambdahat}
 \wh R_\phi(\wh{\lambda}) +2 r\|\wh
\lambda\|_{\ell_1}\le\wh R_\phi(\lambda) + 2r\| \lambda\|_{\ell_1}
\end{equation}
applied to the vector of zeros $\lambda=(0,\ldots,0)^\prime$ implies that
$ \|\wh{\lambda}\|_{\ell_1}\le\phi(0)/(2r)=1/(2r).$ This means that we
can restrict our analysis to the set
\[
\Lambda=\{\lambda\in\RR^M \dvtx \|\lambda\|_{\ell_1} \le1/(2r)\}.
\]
The aim of this section is to show that $\wh{\lambda}$ is close to
$\lambda(r)$ for a judiciously chosen tuning parameter~$r$.

\begin{theorem} \label{cons}
If, for some $p\ge1$,
%
%
\begin{equation}\label{r}
r \ge\frac{1-d}{d} C_F \Biggl\{ 9 \sqrt{ \frac{2\log2(M\vee n)}{n}} + 2
\frac{ p\log_2 n}{\sqrt{2M\vee2n}} + \sqrt{ \frac{2\log 1/\delta}{n}}
\Biggr\},
\end{equation}
then for all $\theta\in\Lambda$, with probability larger than
$1-\delta$,
\[
\Delta R_\phi({\sf f}_{\wh{\lambda}}) + r \|\wh{\lambda}\|_{\ell_1}
\le
\Delta R_\phi({\sf f}_\theta) + 3 r \|\theta\|_{\ell_1} +n^{-p}
\]
and, moreover,
\[
\Delta R_\phi({\sf f}_{\wh{\lambda}}) + r\|\wh{\lambda}-\theta\|
_{\ell_1}
\le\Delta R_\phi({\sf f}_\theta) + 4 r \|\theta\|_{\ell_1} +n^{-p}.
\]
\end{theorem}
\begin{pf}
Write $\wh\delta=\wh{\lambda}-\theta$. Let $\eps=r^{-1}n^{-p}$ and define
%
%
\begin{equation}
\wh r = \sup_{\lambda\in\Lambda} \frac{ \{ \wh R_\phi ({\sf
f}_\lambda) - R_\phi({\sf f}_\lambda) \} - \{ \wh R_\phi({\sf
f}_\theta) - R_\phi({\sf f}_\theta) \} }{\|\lambda-\theta\|_{\ell_1} +
\eps}.
\end{equation}
By Propositions \ref{lemma:MCD} and \ref{lb} in Appendix \ref{appendixB},
\[
\PP\{\wh r\le r\} \ge1-\delta
\]
for the choice $r$ given in (\ref{r}). Rewriting the inequality
(\ref{lambdahat}), we find that
%
%
\begin{eqnarray}\label{WAT}
R_\phi({\sf f}_{\wh{\lambda}}) &\le& R_\phi({\sf f}_\theta) + \{ \wh
R_\phi({\sf f}_\theta) - R_\phi({\sf f}_\theta)  \} - \{ \wh
R_\phi({\sf f}_{\wh{\lambda}}) - R_\phi({\sf f}_{\wh \lambda})
\}\nonumber\\
 &&{}  + 2r\|\theta\|_{\ell_1} - 2r\|\wh{\lambda}\|_{\ell
_1}\\
&\le& R_\phi({\sf f}_\theta) +\wh r (\|\wh\delta\|_{\ell_1} +\eps) +
2r\|\theta\|_{\ell_1} - 2r\|\wh{\lambda}\|_{\ell_1}.\nonumber
\end{eqnarray}
Thus, on the event $\wh r\le r$, after adding
$r\|\wh{\lambda}\|_{\ell_1}$ to both sides, we obtain
\[
R_\phi({\sf f}_{\wh{\lambda}}) + r\|\wh{\lambda}\|_{\ell_1} \le
R_\phi({\sf f}_\theta) +3r\|\theta\|_{\ell_1}+ r\eps,
\]
which proves the first claim. Adding $r\|\wh\delta\|_{\ell_1}$ to both
sides easily yields the second claim.
\end{pf}

A direct consequence of Theorem \ref{cons} is the following corollary
which states that in the sparse setting where
$r\|\lambda(r)\|_{\ell_1}\to0$, the estimator $\wh{\lambda}(2r)$
behaves like the penalized minimizer $\lambda(r)$ in terms of their
risk.

\begin{corollary}\label{cor1}
Suppose that $r\|\lambda(r)\|_{\ell_1}\to0$ as $n\to\infty$ for $r$
satisfying (\ref{r}). Then, with probability at least $1-\delta$,
\[
| \{ R_\phi(\wh{\lambda}) + r\|\wh{\lambda}\|_{\ell_1} \} - \{
R_\phi(\lambda(r))+ r\|\lambda(r)\|_{\ell_1}\}  |\to0
\]
as $n\to\infty$. In particular, when taking $\theta=\lambda(0)$, we
have $ | R_\phi(\wh{\lambda}) - R_\phi(\lambda(0)) |\to0$ and $
\|\wh{\lambda}(2r)-\lambda(0)\|_{\ell_1}= \mathrm{o}(1/r)$.
\end{corollary}
\begin{pf}
We combine the basic property (\ref{lambdahat}) applied to
$\theta=\lambda(r)$ and Theorem \ref{cons}, and we find that on the
event $\wh r\le r$,
\[
R_\phi(\lambda(r)) + r\|\lambda(r)\|_{\ell_1} \le R_\phi(\wh
\lambda)
+ r\|\wh{\lambda}\|_{\ell_1}\le R_\phi(\lambda(r)) + r\|\lambda(r)\|_{\ell_1} + \{
2r\|\lambda(r)\|_{\ell_1} + r\eps\}.
\]
The result then follows from $\{ 2r\|\lambda(r)\|_{\ell_1} + r\eps\}
\to
0$.
\end{pf}



We emphasize that the above results do not impose any restrictions on
the dictionary $\{f_j\}$. If we are willing to make assumptions on the
Gram matrix $\Psi$, then we obtain a more refined result.

\begin{theorem}\label{empirical}
For all $r$ satisfying (\ref{r}) and $\theta\in\Lambda$ such that
$\kappa=\kappa(\theta,7)>0$ and
%
%
\begin{equation}\label{reg2}
\bigl(\kappa^2 r^{\alpha/(1+\alpha)} \| \theta\|_{\ell_0}  \bigr)^{
(1+\alpha)/(2+\alpha) } < c
\end{equation}
for some (small) $c$ depending on $C_F$, $\alpha$, $A$ and $d$, we
have, for some $C$ depending on $c$, that
\[
\Delta R_\phi({\sf f}_{\wh{\lambda}}) + \tfrac12 r
\|\wh{\lambda}-\theta\|_{\ell_1} \le 3\Delta R_\phi(\theta)+
C \| {\sf f}_{\theta}-f_0\|_{\infty}
(\kappa^{-2} r^2 \|\theta\|_{\ell_0})^{(1+\alpha)/(2+\alpha)} +n^{-p}
\]
holds with probability at least $1-\delta$.
\end{theorem}
\begin{pf}
Recall that $\eps= r^{-1}n^{-p}$. We may assume without loss of
generality that
%
%
\begin{equation}\label{een}
R_\phi({\sf f}_\theta)+\eps r \le R_\phi({\sf f}_{\wh{\lambda}
})+\tfrac12
r\|\wh\delta\|_{\ell_1}
\end{equation}
holds, since otherwise the statement holds trivially. Consequently, on
the event $\wh r\le r$, using (\ref{WAT}) and (\ref{een}), we get
\begin{eqnarray*}
R_\phi({\sf f}_{\wh{\lambda}}) &\le& R_\phi({\sf f}_\theta)+ \eps r
+ r
\|\wh\delta\|_{\ell_1} +2 r\|\theta\|_{\ell_1} - 2r\|\wh{\lambda}\|
_{\ell
_1}\\
&\le& R_\phi({\sf f}_{\wh{\lambda}})+\tfrac32 r\|\wh\delta\|_{\ell
_1} +
2r\|\theta\|_{\ell_1}-2r \|\wh{\lambda}\|_{\ell_1}\\
&=& R_\phi({\sf f}_{\wh{\lambda}})+\tfrac32 r \|\wh\delta_I\|_{\ell
_1} +
\tfrac32 r \|\wh{\lambda}_{I^c}\|_{\ell_1}
+2r \|\theta\|_{\ell_1} -2r\|\wh{\lambda}\|_{\ell_1}\\
&=& R_\phi({\sf f}_{\wh{\lambda}}) + \tfrac32 r \|\wh\delta_I\|
_{\ell_1} +
2r\|\theta\|_{\ell_1} -2 r\|\wh{\lambda}_I\|_{\ell_1} -\tfrac12 r\|
\wh
\delta_{I^c}\|_{\ell_1}\\
&\le& R_\phi({\sf f}_{\wh{\lambda}}) + \tfrac72 r
\|\wh\delta_I\|_{\ell_1} -\tfrac12 r\|\wh{\lambda}_{I^c}\|_{\ell_1}
\end{eqnarray*}
so that $\|\wh\delta_{I^c}\|_{\ell_1}\le7\|\wh\delta_{I}\|_{\ell_1}$,
where $I$ is the support of $\theta$. On the other hand,
\begin{eqnarray*}
R_\phi({\sf f}_{\wh{\lambda}}) +\tfrac12 r \|\wh\delta\|_{\ell_1}
&\le&
R_\phi({\sf f}_\theta) +
\eps r + \tfrac32 r \|\wh\delta\|_{\ell_1} +2 r\|\theta\|_{\ell_1} -
2r\|\wh{\lambda}\|_{\ell_1}\\
&\le& R_\phi({\sf f}_\theta) +
\eps r + \tfrac32 r \|\wh\delta_I\|_{\ell_1} +2 r\|\wh\delta_I\|
_{\ell
_1}-{\tfrac12}r\|\wh{\lambda}_{I^c}\|_{\ell_1}\\
&\le& R_\phi({\sf f}_\theta) + {\tfrac72}r \|\wh\delta_{I}\|_{\ell
_1} +
r\eps.
\end{eqnarray*}
The remainder of the proof follows that of Theorem \ref{theorie}, with
$\kappa=\kappa(\theta,7)$.
\end{pf}

This result differs from \cite{W07} (and \cite{TG06} for the case
$d=1/2$) in the appearance of the norm $\| f_0 - {\sf f}_\theta\|
_{\infty}$ on the right-hand side of the (oracle) inequality. This
implies that for $f_0={\sf f}_\theta$ and for some sparse $\theta
=\lambda(0)$ satisfying the conditions of Theorem \ref{empirical}, we
can expect fast rates, regardless of the classification complexity!
Another important difference with both papers is that no restriction is
imposed on the sup-norm of ${\sf f}_\lambda$. Such a condition is
unnatural as $|{\sf f}_\lambda|\le C$ may overrule the restriction that
the penalty term $r\|\lambda\|_{\ell_1}$ imposes. \\

We now consider bounds on the error and reject rates without an
additional test sample. We write
\[
\PP_n\{Y {\sf f}_{\wh{\lambda}}(X) \le\beta\} = \frac1n \sum_{i=1}^n
I\{ Y_i {\sf f}_{\wh{\lambda}}(X_i) \le\beta\}
\]
for any $\beta>0$. The misclassification and rejection rate can be
bounded above as follows.

\begin{theorem}\label{theorem:errorrate}
If
\[
r(\gamma) \ge\frac{9 C_F}{\gamma} \sqrt{\frac{2\log2(M\vee n)}{n}}
+\frac{2p \log_2(n) C_F }{\gamma\sqrt{2(M\vee n)}} + \frac
{C_F}{\gamma
} \sqrt{ \frac{ 2\log(1/\delta)}{n}},
\]
then, with probability at least $1-\delta$, we have
\begin{eqnarray*}
\PP\{ Y {\sf f}_{\wh{\lambda}}(X) \le-\tau\} &\le& \min_{\gamma>0} [
\PP_n\{Y {\sf f}_{\wh{\lambda}}(X) \le-\tau+ \gamma\} +r(\gamma)
\|\wh{\lambda}\|_{\ell_1}
] + n^{-p},\\
\PP\{ | {\sf f}_{\wh{\lambda}}(X) |\le\tau\} &\le& \min_{\gamma>0} [
\PP_n\{| {\sf f}_{\wh{\lambda}}(X)| \le\tau+ \gamma\} +r(\gamma)
\|\wh{\lambda}\|_{\ell_1}] + n^{-p}.
\end{eqnarray*}
\end{theorem}
\begin{pf}
Set
\[
\varphi_\gamma(z) =
\cases{
1, &\quad if $z<-\tau$,\vspace*{1pt}\cr
\dfrac{1}{\gamma} (\gamma-\tau-z), &\quad if $-\tau\le z \le-\tau+\gamma$,\vspace*{2pt}\cr
0, &\quad if $z\ge-\tau+\gamma$.
}
\]
The following inequalities then hold uniformly in $\lambda$:
\begin{eqnarray*}
\PP\{ Y {\sf f}_ \lambda(X) \le-\tau\} &\le& \PP_n\{Y {\sf
f}_\lambda
(X) \le-\tau+ \gamma\} +
R_{\varphi_\gamma}({\sf f}_\lambda) - \wh R_{\varphi_\gamma}({\sf
f}_\lambda)\\
&\le& \PP_n\{Y {\sf f}_\lambda(X) \le-\tau+ \gamma\} +\wh r_0
\{ \|\lambda\|_{\ell_1}+\eps\},
\end{eqnarray*}
where
\[
\wh r_0= \sup_{\lambda\in\Lambda} \frac{ | \wh R_{\varphi _\gamma}
({\sf f}_\lambda) - R_{\varphi_\gamma} ({\sf f}_\lambda) | }{ \|
\lambda\|_{\ell_1} + \eps},
\]
with $\eps$ given by $\eps r(\gamma)= n^{-p}$. We can invoke
Propositions \ref{lemma:MCD} and \ref{lb} to complete the proof of the
first claim. The proof of the second claim uses the reasoning above,
with the only modification being that $\varphi_\gamma(z)$ is now given
by
\[
\varphi_\gamma(z) =
\cases{
1, &\quad if $|z|<\tau$,\vspace*{1pt}\cr
\dfrac{1}{\gamma} (z+\gamma+\tau), &\quad if $-\tau-\gamma\le z
\le
-\tau$,\vspace*{1pt}\cr
-\dfrac{1}{\gamma} (z-\gamma-\tau), &\quad if $\tau\le z \le\tau
+\gamma$,\vspace*{2pt}\cr
0, &\quad if $|z|\ge\tau+\gamma$;
}
\]
the rest of the reasoning is unchanged.
\end{pf}
%

\section{Numerical experiments}\label{sec4}

\begin{figure}

\includegraphics{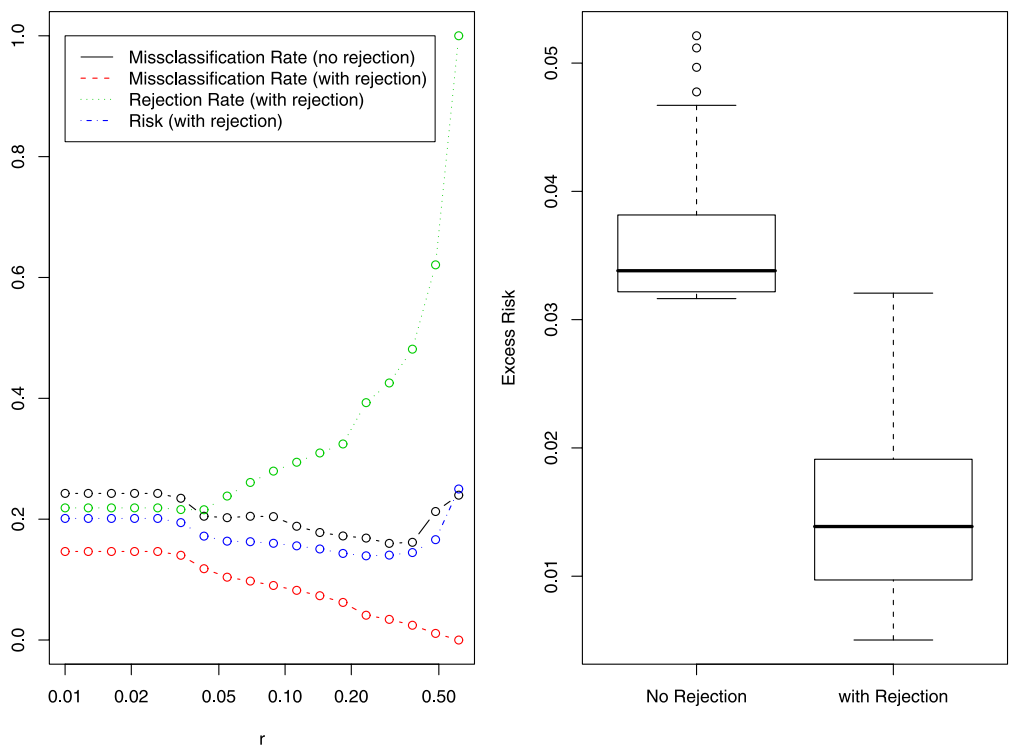}

\caption{Simulation -- the effect of rejection, misclassification rate
and excess risk $R_\ell$. The left-hand panel shows the three criteria
as functions of the tuning parameter $r$ for the support vector machine
(SVM) with rejection option for a typical run. Also included is the
misclassification rate for the usual SVM. It is evident that SVM with
rejection option enjoys lower misclassification rate by withholding
decision for ``hard-to-classify'' cases. The right-hand panel compares
the excess $\ell$-risk for SVM with or without rejection option. The
box plots of the excess risk are produced based on 200 runs. This again
confirms that SVM with rejection option leads to improved performance
in terms of the $\ell$ loss.} \label{fig:sim}
\end{figure}

We now demonstrate the practical merits of $\wh{\lambda}(r)$ via a
couple of numerical experiments. We begin by noting that the
computation of $\wh{\lambda}(r)$ can be conveniently formulated as a
linear program. Let $\xi_1,\ldots,\xi_n$ be the slack variables such
that
%
%
\begin{equation}
\xi_i\ge0,\qquad  \xi_i\ge1-Y_if_\lambda(X_i),\qquad \xi_i\ge 1-aY_if_\lambda(X_i).
\end{equation}
Clearly the minimum $\xi_i$ that satisfies these constraints is
$\phi(Y_if(X_i))$. We also introduce slack variables $\xi_{n+i}$,
$i=1,\ldots, M,$ to represent $|\lambda_i|$, that is,
%
%
\begin{equation}
\xi_{n+i}\ge\lambda_i,\qquad \xi_{n+i}\ge-\lambda_i.
\end{equation}
Using the slack variables, $\wh{\lambda}(r)$ can be given as the solution
of the linear program
\begin{eqnarray*}
\min_{\lambda, \xi} [\xi_1+\cdots+\xi_n+r(\xi_{n+1}+\cdots+\xi _{n+M})]
\end{eqnarray*}
subject to
\begin{eqnarray*}
\xi_i&\ge&0,\qquad \xi_i\ge1-y_ih_i,\qquad \xi_i\ge1-aY_ih_i,\qquad   i=1,\ldots,
n,\\
\xi_{n+i}&\ge&\lambda_i,\qquad \xi_{n+i}\ge-\lambda_i,\qquad   i=1,\ldots
, M,\\
h_i&=&\sum_j \lambda_jf_j(X_i),\qquad   i=1,\ldots, n.
\end{eqnarray*}
To illustrate the merits of $\wh{\lambda}$, we implement the method
described above and first apply it to a set of simulated examples. To
fix ideas, we set $d=0.25$ or, equivalently, $a=3$. For each run, $50$
positive instances ($Y=+1$) and $50$ negative instances ($Y=-1$) were
generated. Two hundred ($M=200$) features ($f_j$'s) were simulated from
a multivariate normal distribution. For positive instances, the mean
was set to $(1/\sqrt{2}, 1/\sqrt{2},0,\ldots,0)',$ whereas for the
negative instances, the mean was set to
$(-{1/\sqrt{2}},-1/\sqrt{2},0,\ldots,0)'$. In both cases, the
covariance matrix was the identity matrix. The operating
characteristics of the method are demonstrated in Figure \ref{fig:sim}.
On the left-hand side, the misclassification rate (${\mathbb P}(Y{\sf
f}_{\hat{\lambda}}(X)<-0.5)$), rejection rate (${\mathbb P}(|Y{\sf
f}_{\hat{\lambda}}(X)|<0.5)$) and associated $\ell$-risk of the
$\ell_1$-regularized generalized hinge loss ($R_\ell({\sf
f}_{\hat{\lambda}})$) are plotted as functions of the tuning parameter
$r$ for a typical simulation. The results are to be compared with the
usual $\ell_1$-regularized support vector machines where no rejection
option is allowed. Since there is no rejection, the misclassification
rate for the usual support vector machines coincides with its
$\ell$-risk. It is evident that by incorporating the rejection option,
$\wh{\lambda}$ yields a smaller $\ell$-risk, provided that both methods
are optimally tuned. To further investigate the merits of allowing the
rejection option, we repeated the experiment 200 times. The excess
risk $\Delta R_\ell$ of both the usual support vector machine and the
proposed method are summarized in the plot on the right-hand side. It
further confirms the advantage of $\wh{\lambda}$.

%

\begin{figure}

\includegraphics{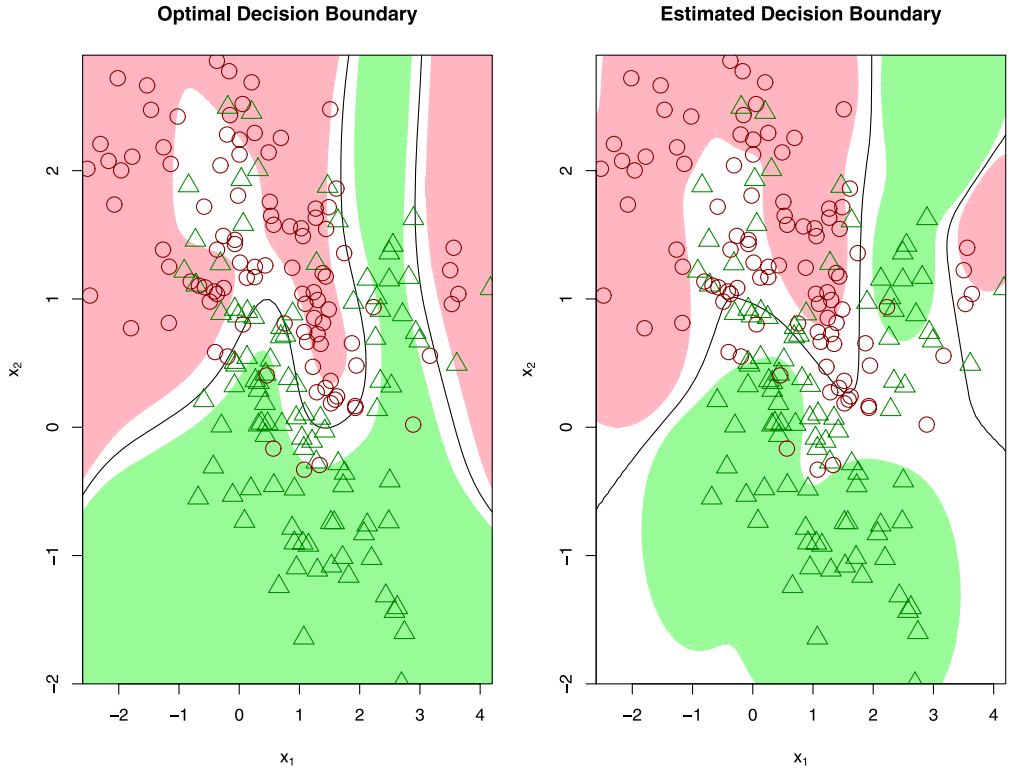}

\caption{Mixture data -- optimal and estimated decision boundaries. The
left-hand panel gives the optimal decision boundary, whereas the
right-hand panel corresponds to the SVM with rejection option. In both
plots, positive cases are represented by red circles and negative cases
by green triangles. The light red regions correspond to classification
$Y=+1$ and light green regions to classification $Y=-1$. Areas where a
decision is withheld are not shaded. The solid black line in the
left-hand panel is the level set for $\PP(Y=+1|X)=0.5$. The solid black
line in the right-hand panel is the level set for ${\sf
f}_{\wh{\lambda}}=0$.} \label{fig:mix}
\end{figure}

To further demonstrate the merits of the method, we apply it to the
mixture data example considered in \cite{HTF}. The training data
consist of 200 data points generated from a pair of two-dimensional
mixture densities. Similarly to \cite{HTF}, we consider a dictionary of
Gaussian radial basis functions $ f_j(\cdot)=\exp(-2\|\cdot-b_j\|^2),
 j=1,\ldots, 100, $ where the locations $b_j$ are placed on a~$10\times10$ equally spaced lattice. To fix ideas, we consider the
case where $d=0.25$. The optimal classification rule will classify an
observation as $+1$ if the corresponding conditional probability
$\PP(Y=+1|X)$ is greater than $0.75$ and as $-1$ if the conditional
probability is less than $0.25$. When the conditional probability is
between $0.25$ and $0.75$, we withhold the decision. The corresponding
decision boundaries are given in the right-hand panel of Figure
\ref{fig:mix}. It is known that the usual SVM only targets the decision
boundary identified with $\PP(Y=+1|X)$ and cannot be used to recover
the optimal decision boundaries given here; see, for instance,
\cite{YW09} for further discussion of this issue. In contrast, the SVM
with rejection option is devised specifically for this purpose. To this
end, we ran the SVM with rejection option with $a=3$ and $\tau=0.5,$ as
discussed earlier. The tuning parameter $r$ was selected by tenfold
cross-validation. The left-hand panel of Figure \ref{fig:mix} gives the
estimated decision boundaries. It is clear from the plot that SVM with
rejection option successfully captured the main characteristics of the
underlying probabilities. The main difference between the two sets of
decision boundaries occurs in regions where no observations are
available. As a result, the SVM with rejection option opted for
withholding a decision.

%

\begin{appendix}\label{appendix}
\section{Connection between excess risk and weighted $L_2$~norm}\label{appendixA}

The next lemma is a technical result that links the excess risk $\Delta
R_\phi(\lambda)$ to the $L_2$ norm:
\[
\| {\sf f}_\lambda- f_0\|=\sqrt{ \EE[ | {\sf f}_\lambda(X)-f_0(X)
|^2 \omega(X) ] }
\]
with
$\omega(X)= \eta(X)(1-\eta)(X)$. Its proof is rather technical and
relies on results obtained in \cite{BW08}. Essentially, $\| {\sf
f}_\lambda- f_0\|_{\infty}$ replaces the suboptimal bound $1+
C_\Lambda C_F$ in \cite{W07}.

\begin{lemma}\label{T}
Let $\alpha>0$ be as in Definition 2.2. Then,
for all $\lambda\in\RR^M$,
%
%
\begin{equation}
\| {\sf f}_\lambda- f_0 \| ^{ 2+2\alpha} \le4A (2d)^\alpha\| {\sf
f}_\lambda- f_0\|_{\infty} ^{2+\alpha} \{ \Delta R_\phi(\lambda) \}
^\alpha.
\end{equation}
\end{lemma}

\begin{pf} Let $f\dvtx\mathcal X\to\RR$ be arbitrary and
set
\[
\rho_\eta(f,f_0) =
\cases{
\eta|f-f_0|, &\quad if $\eta<d$ and $f<-1$, \cr
(1-\eta)|f-f_0|, &\quad if $\eta>1-d$ and $f>1$, \cr
|f-f_0|, &\quad otherwise,
}\label{drie}
\]
then \cite{BW08}, Lemma 9, states that
\[
\Delta R_\phi(\lambda) \ge d^{-1} \EE\bigl[ \rho_\eta(f,f_0)(X)
\bigl(|\eta(X)-(1-d)| I_{\{X\in E_-\}} + |\eta(X)-d| I_{\{X\in E_+\}}\bigr) \bigr]
\]
with
\[
E_-=\{|\eta-(1-d)|\le|\eta-d|\},\qquad  E_+=\{|\eta-(1-d)|>|\eta -d|\}.\
\]
Using (A.1), for any set $E$,
\begin{eqnarray*}
&& \EE\bigl[ \rho_\eta(f,f_0)(X) |\eta(X) -(1-d)| I_{\{X\in E \}} \bigr]\\
&&\quad\ge t \EE\bigl[ \rho_\eta(f,f_0)(X) I_{\{ |\eta(X)-(1-d)|\ge t\}} I_{\{
X\in
E\}} \bigr] \\
&&\quad= t \EE\bigl[ \rho_\eta(f,f_0)(X) I_{\{ X\in E\}} \bigr] - t \EE\bigl[
\rho_\eta(f,f_0)(X) I_{\{ |\eta(X)-(1-d)| < t, X\in E\}} \bigr]\\
&&\quad\ge t \EE\bigl[ \rho_\eta(f,f_0)(X) I_{\{ X\in E\}} - \|f-f_0\|_{\infty} A
t^\alpha\bigr].
\end{eqnarray*}
Similarly,
\[
\EE\bigl[ \rho_\eta(f,f_0)(X) |\eta(X)-d| I_{\{ X\in E\}} \bigr] \ge t \EE\bigl[
\rho_\eta(f,f_0)(X) I_{\{ X\in E\}} - \|f-f_0\|_{\infty} A t^\alpha\bigr],
\]
and we obtain
\begin{eqnarray*}
\Delta R_\phi(\lambda) &\ge& d^{-1} t \EE\bigl[ \rho_\eta({\sf
f}_\lambda,f_0)(X) I_{\{ X\in E_+\cup E_-\} }- 2\|{\sf
f}_\lambda-f_0\|_{\infty} At^\alpha
\bigr]\\
&=& d^{-1} t \EE[ \rho_\eta({\sf f}_\lambda,f_0)(X) -2\|{\sf
f}_\lambda-f_0\|_{\infty} At^\alpha].
\end{eqnarray*}
Plugging
\[
t = \biggl( \frac{ \EE[ \rho_\eta({\sf f}_\lambda,f_0)(X)] }{4 A \| {\sf
f}_\lambda-f_0\|_{\infty} } \biggr)^{1/\alpha}
\]
into the
preceding expression, we obtain
\[
\Delta R_\phi(\lambda) \ge\frac{( \EE[ \rho_\eta ({\sf
f}_\lambda,f_0)(X) ] )^{(1+\alpha)/\alpha} } {2d(4A\| {\sf
f}_\lambda-f_0\|_{\infty} )^{1/\alpha}} .
\]
Since
\[
\| {\sf f}_\lambda- f_0\|^2 = \EE[ \omega(X) ({\sf f}_\lambda
-f_0)^2(X) ] \le\| {\sf f}_\lambda-f_0\|_{\infty}\EE[ \omega(X) |{\sf
f}_\lambda(X)-f_0(X)|],
\]
we get,
for all $\lambda$,
\begin{eqnarray*}
\Delta R_\phi(\lambda) &\ge& \frac{(\EE[ \omega(X) | {\sf f}_\lambda(X)
-f_0(X) |  ])^{(1+\alpha)/\alpha}}{2d(4A \|{\sf
f}_\lambda-f_0\|_{\infty}
)^{1/\alpha}} \\
&\ge& \frac{(\| {\sf f}_\lambda-f_0\| ^2)^{(1+\alpha)/\alpha}}{2d(4A
)^{1/\alpha} \|{\sf
f}_\lambda-f_0\|_{\infty}^{(2+\alpha)/\alpha}} .
\end{eqnarray*}
The claim follows.
\end{pf}

\begin{remark}
If $|{\sf f}_\lambda|\le1$, then $\rho_\eta({\sf f}_\lambda
,f_0)=|{\sf
f}_\lambda- f_0|$. Hence, if we restrict the parameters $\lambda$ such
that ${\sf f}_\lambda$ are bounded by 1, then we can impose the
restricted eigenvalue condition on the matrix with entries $\EE[ f_i(X)
f_j(X)]$ instead of $\EE[ f_i(X) f_j(X) \omega(X)]$.
\end{remark}

\section{A maximal inequality for a weighted empirical~process}\label{appendixB}

Recall that $\Lambda=\{ \lambda\in\RR^M\dvtx \|\lambda\|_{\ell _1}\le 1/(2
r) \}$ and let $\theta\in\Lambda$ and $\eps>0$. Let $\varphi\dvtx\RR\to\RR$
be a convex function with Lipschitz constant $C_\varphi$ and define the
risks
\begin{eqnarray*}
R_\varphi({\sf f}_\lambda) &=& \EE[ \varphi(Y {\sf f}_\lambda
(X))],\\
\wh R_\varphi({\sf f}_\lambda) &=& \frac1n \sum_{i=1}^n \varphi(Y_i
{\sf f}_\lambda(X_i)).
\end{eqnarray*}

Finally, let $\eps>0$ and set
\[
\wh r(\varphi,\theta,\eps) =\sup_{\lambda\in\Lambda} \frac{
|\{ \wh R_\varphi({\sf f}_\lambda)- R_\varphi({\sf f}_\lambda) \} -\{
\wh R_\varphi({\sf f}_\theta) - R_\varphi({\sf f}_\theta)\} |}{
\|\theta-\lambda\|_{\ell_1} +\eps}.
\]
We prove a maximal inequality for $\wh r(\varphi,\theta,\eps)$ which
slightly generalizes the result obtained in~\cite{W07}.

\begin{prop} \label{lemma:MCD}
Let $0<\delta<1$ and set
\[
r(\varphi,\theta,\eps)= \EE[ \wh r(\varphi,\theta,\eps) ] +C_\varphi
C_F \sqrt{ \frac{2\log (1/\delta)}{n}}.
\]
Then,
\[
\PP\{ r(\varphi,\theta,\eps)\ge\wh
r(\varphi,\theta,\eps)\}\ge1-\delta.
\]
\end{prop}

\begin{pf}
First, observe that changing a pair $(X_i,Y_i)$ in $\wh r$ changes it
by at most $2C_\varphi C_F/n$. The result follows immediately after
applying McDiarmid's exponential inequality \cite{DL00}, Theorem 2.2, page
8.
\end{pf}

We now control the expectation of $\wh r(\varphi,\theta,\eps)$.

\begin{prop}\label{lb}
Set $J=\lceil\log_2(1/\{\eps r\}) \rceil$. Then,
\[
\EE[ \wh r(\varphi,\theta,\eps) ] \le 9 C_\varphi C_F
\sqrt{\frac{2\log2(M\vee n)}{n}} +\frac{2J C_\varphi C_F
}{\sqrt{2(M\vee n)}}.
\]
\end{prop}

\begin{pf} 
Let $\sigma_1,\ldots,\sigma_n$ be independent Rademacher variables,
taking the values $\pm1$, each with probability 1$/$2, independent of the
data $(X_1,Y_1),\ldots,(X_n,Y_n)$. Set
\[
\wh R_\varphi^0({\sf f}_\lambda) = \frac1n \sum_{i=1}^n \sigma_i
\varphi(Y_i {\sf f}_\lambda(X_i)).
\]

A standard symmetrization trick \cite{DL00}, page 18, shows that
\begin{eqnarray*}
\EE[ \wh r (\varphi,\theta, \eps) ] &\le&2 \EE \biggl[ \sup
_{\lambda\in\Lambda} \frac{ |\wh R_\varphi^0({\sf f}_\lambda)- \wh
{R}_\varphi^0({\sf f}_\theta) |}{\|\lambda-\theta\|_{\ell_1} + \eps}
\biggr]\\
&\le& 2\EE\biggl[ \sup_{\|\lambda-\theta\|_{\ell_1}\le\eps} \frac{ |\wh
R_\varphi^0({\sf f}_\lambda)- \wh{R}_\varphi^0({\sf f}_\theta)
|}{\|\lambda-\theta\|_{\ell_1} + \eps} \biggr]+2
\EE\biggl[\sup_{\eps\le\|\lambda-\theta\|_{\ell_1}\le1/r} \frac{ |\wh
R_\varphi^0({\sf f}_\lambda)- \wh{R}_\varphi^0(\theta ) |}{\|
\lambda-\theta\|_{\ell_1} + \eps}\biggr]
\\
&=& (I) + (\mathit{II}).
\end{eqnarray*}
The first term
\begin{eqnarray*}
I &=&2 \EE\biggl[ \sup_{\|\lambda-\theta\|_{\ell_1}\le\eps} \frac{ |\wh
R_\varphi^0({\sf f}_\lambda)- \wh{R}_\varphi^0({\sf f}_\theta)
|}{\|\lambda-\theta\|_{\ell_1} + \eps}
\biggr]\\
&=& \EE\Biggl[ \sup_{\|\lambda-\theta\|_{\ell_1}\le\eps} \frac{ 1}{\|
\lambda-\theta\|_{\ell_1} + \eps} \Biggl| \frac1n \sum_{i=1}^n \{ \varphi(
Y_i {\sf f}_\lambda(X_i) ) - \varphi(Y_i {\sf f}_{\theta }(X_i))\} \Biggr| \Biggr]
\end{eqnarray*}
can be bounded using the contraction principle
for Rademacher processes; see \cite{LT}, pages 112--113. For this, we
observe that the function $g(z)= \varphi(z_0+ z)-\varphi(z_0)$ is
Lipschitz with Lipschitz constant $C_\varphi$ and $g(0)=0$.
Consequently,
\begin{eqnarray*}
(I) &\le& 2\frac{C_\varphi}{\eps} \EE\Biggl[ \sup_{\|\lambda -\theta\|
_{\ell_1} \le\eps} \Biggl| \frac1n \sum_{i=1}^n \sigma_i Y_i {\sf
f}_{\lambda- \theta}(X_i) \Biggr| \Biggr]
\\
&\le& 2\frac{C_\varphi}{\eps} \EE\Biggl[ \sup_{\|\lambda-\theta\|_{\ell_1}
\le\eps} \|\lambda -\theta \|_{\ell_1} \max_{1\le j\le M} \Biggl| \frac1n
\sum_{i=1}^n \sigma _i Y_i f_j(X_i) \Biggr| \Biggr]
\\ &\le&2 C_\varphi\EE\Biggl[\max_{1\le j\le M} \Biggl| \frac1n
\sum
_{i=1}^n \sigma_i Y_i f_j(X_i) \Biggr| \Biggr]\\
&\le&2 C_\varphi C_F \frac{\sqrt{2\log(2M)}}{\sqrt{n}}.
\end{eqnarray*}
The last maximal inequality can be found in
\cite{DL00}, Lemma 2.2, page 7, which uses the fact that the variables
$\sigma_i Y_i f_j(X_i)$ are sub-Gaussian,
\[
\EE\Biggl[\exp\Biggl\{ s \sum_{i=1}^n \sigma_i Y_i f_j(X_i) \Biggr\}\Biggr] \le\exp(n s^2 C_F^2
/2)
\]
for all $s$,
which follows, in turn, from \cite{DL00}, Lemma 2.1, page 5.

The second term (II) requires a peeling argument \cite{sara}, page 70.
Since $0\le\wh r\le2 C_\varphi C_F$ almost surely, we can use the bound
%
%
\begin{equation}\label{ong}
(\mathit{II}) \le \zeta+2 C_\varphi C_F \PP\biggl\{ \sup_{\eps\le\| \lambda
-\theta\|_{\ell_1}\le1/r} 2 \frac{| \wh R_\varphi^0({\sf f}_\lambda)
-\wh R_\varphi^0({\sf f}_\theta) |}{\|\lambda- \theta\|_{\ell_1} + \eps
} \ge\zeta\biggr\}.\
\end{equation}
Observe that for any $\zeta>0$,
\begin{eqnarray*}
&& \PP\biggl\{ \sup_{\eps\le\|\lambda-\theta\|_{\ell_1}\le1/r}2 \frac{| \wh
R_\varphi^0({\sf f}_\lambda) -\wh R_\varphi^0({\sf f}_\theta)
|}{\|\lambda- \theta\|_{\ell_1} + \eps} \ge\zeta \biggr\}
\\
&&\quad\le \sum_{j=1}^{J} \PP\biggl\{ \sup_{2^{j-1}
\eps\le\|\lambda-\theta\|_{\ell_1}\le2^j \eps} | \wh R_\varphi^0({\sf
f}_\lambda) -\wh R_\varphi^0({\sf f}_\theta )| \ge2^{j-2} \eps\zeta\biggr\}.
\end{eqnarray*}
Now, set
\[
Z_j= \sup_{ \|\lambda- \theta\|_{\ell_1}\le2^j \eps} | \wh
R_\varphi^0({\sf f}_\lambda) -\wh R_\varphi^0({\sf f}_ \theta )|
\]
and the same considerations leading to the final bound of (I) above
yield
\[
\EE[ Z_j ] \le2^j \eps C_\phi C_F \frac{ \sqrt{2 \log
(2M)}}{\sqrt{n}}
\]
and for $t=1/\sqrt{2}$, we obtain
\[
(\mathit{II}) \le \zeta+
2C_\varphi C_F \sum_{j=1}^{J} \PP\{ Z_j -\EE[ Z_j ] \ge 2^{j-2}
\eps\zeta-\EE[ Z_j ]\}.
\]
A change of a single pair $(X_i,Y_i)$ changes $Z_j$ by at most
$2C_\varphi C_F (2^j \eps) /n$, so that another application of the
bounded differences inequality \cite{DL00}, Theorem 2.2, page 8, gives,
by taking
\[
\zeta= 7 C_\varphi C_F \frac{ \sqrt{2\log2(M\vee n)} }{\sqrt{n}},
\]
the final bound
\begin{eqnarray*}
&& \sum_{j=1}^{J} \PP\{ Z_j -\EE[ Z_j ] \ge2^{j-2} \eps\zeta -\EE
[ Z_j ]\}\\
&&\quad\le \sum_{j=1}^{J} \PP\biggl\{ Z_j -\EE[ Z_j ] \ge t\cdot2^{j} C_\varphi
C_F \eps\frac{\sqrt{2\log(2M\vee 2n)}}{\sqrt{n}}
\biggr\}\\
&&\quad\le J \exp\biggl\{-2 \frac{t^2 ( C_\varphi C_F 2^j \eps)^2 2\log
(2M\vee2n)}{(2C_\phi C_F 2^j \eps)^2 }\biggr\}\\
&&\quad= J (2 M\vee2n)^{-t^2} < J /\sqrt{ 2M\vee2n}.
\end{eqnarray*}
Finally, we invoke (\ref{ong}) to complete the proof of Proposition
\ref{lemma:MCD}.
\end{pf}
\end{appendix}

\section*{Acknowledgements} The research of Marten Wegkamp was
supported in part by NSF Grant DMS-0706829. The research of Ming Yuan
was supported in part by NSF Grant DMS-08-46234 and NIH Grant R01GM076274-01.

\printhistory


\begin{thebibliography}{17}

\bibitem{BW08}
Bartlett, P.L. and Wegkamp, M.H. (2008). Classification with a reject
option using a hinge loss. \textit{J.~Mach. Learn. Res.}
\textbf{9} 1823--1840.
\MR{2438825}

\bibitem{BRT}
Bickel, P.J., Ritov, Y. and Tsybakov, A.B. (2009). Simultaneous
analysis of Lasso and Dantzig selector. \textit{Ann. Statist.} \textbf{37}
1705--1732.
\MR{2533469}

%
%

\bibitem{DL00}
Devroye, L. and Lugosi, G. (2000). \textit{Combinatorial Methods in
Density Estimation}. New York: Springer.
\MR{1843146}

\bibitem{HTF}
Hastie, T., Tibshirani, R. and Friedman, J. (2001). \textit{The
Elements of
Statistical Learning.} New York: Springer.
\MR{1851606}

\bibitem{HW06}
Herbei, R. and Wegkamp, M.H. (2006). Classification
with reject option. \textit{Canad. J. Statist.} \textbf{34}
709--721.
\MR{2347054}

\bibitem{K09}
Koltchinskii, V. (2009). Sparsity in penalized empirical risk
minimization. \textit{Ann. Inst. H.
Poincar\'{e} Probab. Statist.} \textbf{45} 7--57.
\MR{2500227}

\bibitem{LT}
Ledoux, M. and Talagrand, M. (1991). \textit{Probability in
Banach Spaces.} New York: Springer.
\MR{1102015}


\bibitem{TG06}
Tarigan, B. and van de Geer, S.A. (2006).
Classifiers of support vector machine type with $\ell_1$ complexity
regularization.
\textit{Bernoulli} \textbf{12} 1045--1076.
\MR{2274857}

\bibitem{T04}
Tsybakov, A.B. (2004). Optimal aggregation of
classifiers in statistical learning. \textit{Ann. Statist.}
\textbf{32} 135--166.
\MR{2051002}

\bibitem{sara}
van de Geer, S.A. (2000). \textit{Empirical Processes in M-estimation.}
Cambridge: Cambridge Univ. Press.


\bibitem{W07}
Wegkamp, M.H. (2007). Lasso type classifiers with a reject
option. \textit{Electron. J. Statist.} \textbf{1} 155--168.
\MR{2312148}

\bibitem{YW09}
Yuan, M. and Wegkamp, M.H. (2010). Classification methods
with reject option based on convex risk minimization. \textit{J.
Mach. Learn. Res.} \textbf{11} 111--130.
\MR{2591623}

\end{thebibliography}
\end{document}